\documentclass[10pt,a4paper]{article}[1.00]
\usepackage{fancyhdr}
\usepackage{textcomp}

\usepackage[english,UKenglish]{babel}
\usepackage[latin1]{inputenc} % Umlaute verwenden
\usepackage{theorem}

\usepackage[T1]{fontenc} % 8 Bit Fonts erleichtern Trennung
\usepackage{amsmath}
\usepackage{amssymb}
\newcommand{\quat}[1]{\mbox{$\left(
    \begin{array}{c} #1_0\\
				     #1_1\\
					 #1_2\\
					 #1_3
				 \end{array}
				\right)$}}

\newcommand{\quatvek}[1]{\mbox{$\left(
    \begin{array}{c} 0\\
				     #1_1\\
					 #1_2\\
					 #1_3
				 \end{array}
				\right)$}}

\newcommand{\vek}[1]{\mbox{$\left(
     \begin{array}{c} #1_1\\
	                  #1_2\\
					  #1_3
				\end{array} 
				\right)$}}	
					
\newcommand{\vekk}[1]{\mbox{$\left(
      \begin{array}{c}#1_0\\
	               \mathbf{#1}
		\end{array}\right)$}}

\newcommand{\vekd}[1]{\mbox{$\left(
      \begin{array}{c}0\\
	               \mathbf{#1}
		\end{array}\right)$}}

\newcommand{\einsq}{\mbox{$\left(
       \begin{array}{c} 1\\
	                    0\\
						0\\
						0
		\end{array}\right)$}}

\newcommand{\einsqv}{\mbox{$\left(
      \begin{array}{c}1\\
	               \mathbf{0}
		\end{array}\right)$}}		
\newcommand{\ixsqv}[1]{\mbox{$\left(
      \begin{array}{c}#1\\
	               \mathbf{0}
		\end{array}\right)$}}

\newcommand{\Ematrix}{\mbox{$\left(
      \begin{array}{*{3}{r@{\;}}r}
	        1 & 0 & 0 & 0\\
			0 & 1 & 0 & 0\\
			0 & 0 & 1 & 0\\
			0 & 0 & 0 & 1     
      \end{array}\right)$}}		
\newcommand{\qmatrix}[1]{\mbox{$\left(
      \begin{array}{*{3}{r@{\;}}r}
	        #1_0 & -#1_1 & -#1_2 & -#1_3\\
			#1_1 &  #1_0 & -#1_3 &  #1_2\\
			#1_2 &  #1_3 &  #1_0 & -#1_1\\
			#1_3 & -#1_2 &  #1_1 &  #1_0     
      \end{array}\right)$}}	
	
\newcommand{\admatrixy}[2]{\mbox{$\left(
      \begin{array}{*{3}{c@{\;}}c}
	        #2-#1_0 & -#1_1 & -#1_2 & -#1_3\\
			#1_1 &  #2-#1_0 & -#1_3 &  #1_2\\
			#1_2 &  #1_3 &  #2-#1_0 & -#1_1\\
			#1_3 & -#1_2 &  #1_1 &  #2-#1_0     
      \end{array}\right)$}}

\newcommand{\quatSang}[1]{\mbox{$#1_0+#1_1$\mathversion{bold}$\mathit{i}$\mathversion{normal}
                                  $+#1_2$\mathversion{bold}$\mathit{j}$\mathversion{normal}
								  $+#1_3$\mathversion{bold}$\mathit{k}$\mathversion{normal}}}		

\begin{document}
\title{\bf N-th root}
\author{Jochen Hans\\ Liliengasse 3\\ 01067 Dresden\\ math@jochenhans.de}
\maketitle
\begin{abstract} The quaternion equation $\mathbf{X}^n=\mathbf{A}$ is solved for any integer number
$n\geq 2$.
$\mathbf{A}$ is a given quaternion with komplex numbers as its elements. 
We use the isomorphism between quaternions and $(4,4)$-matrices to solve this equation.
\end{abstract}
\section{Notations}
latin lower case letter -  complex number\\
$i$ denotes here the imaginary unit ($i^2=-1$)\\
bold latin lower case letter -  a 3-dim. vector,\\ i. e. $\mathbf{a}=\vek{a}$\\
bold latin capital letter  - a quaternion,\\ i. e. $\mathbf{A}=\quat{a}=\vekk{a}$\\
A quaternion is often written as $\quatSang{a}$ but in this article I use the above definition.\\
capital Gothic type -  a matrix   i.e.  $\mathfrak{A}=\qmatrix{a}$\\
The inner product - $\mathbf{a}\mathbf{b}:=a_1b_1+a_2b_2+a_3b_3$.\\
The quadrat of a vector - $\mathbf{a}^2:=a_1^2+a_2^2+a_3^2$\\
The quaternion product -
\begin{equation} \mathbf{A}\mathbf{B}=\quat{a}.\quat{b}=\left(\begin{array}{c}
               a_0b_0-a_1b_1-a_2b_2-a_3b_3\\
			   a_0b_1+a_1b_0+a_2b_3-a_3b_2\\
			   a_0b_2-a_1b_3+a_2b_0+a_3b_1\\
			   a_0b_3+a_1b_2-a_2b_1+a_3b_0
          	     \end{array}\right)\end{equation}
or
\begin{equation}\mathbf{A}\mathbf{B}=
          \vekk{a}\vekk{b}=\left(\begin{array}{c}
		                  a_0b_0-\mathbf{a}\mathbf{b}\\
						  a_0 \mathbf{b}+b_0 \mathbf{a}+\mathbf{a}\times\mathbf{b}
						  \end{array}\right).\end{equation}
The isomorphism between quaternions and $4\times 4$ matrices
$f_I: \mathbb{Q}\longrightarrow\mathbb{M}^{4\times 4}$ is:\\[1.5ex]
$f_I(\quat{a})=\qmatrix{a}$\\[1.5ex]
There are another isomorphisms too but they do not yield another result.\\
The characteristic polynomial of a matrix: $c(y):=\mbox{det}(y\mathfrak{E}-\mathfrak{A})$.\\
The minimal polynomial of a matrix: $m(y):=$ annihilating polynomial of least degree.\\
The adjugate to $(y\mathfrak{E}-\mathfrak{A})$: $\mathfrak{B}(y):=\mbox{adj}(y\mathfrak{E}-\mathfrak{A})$.\\
This relation is valid - $(y\mathfrak{E}-\mathfrak{A})\mathfrak{B}(y)=c(y)\mathfrak{E}$\\
The reduced adjugate to $(y\mathfrak{E}-\mathfrak{A})$: $\mathfrak{C}(y)$ with
$(y\mathfrak{E}-\mathfrak{A})\mathfrak{C}(y)=m(y)\mathfrak{E}$
\section{Method of solution}
We have to solve the quaternion equation
\begin{equation}\label{xhochnq}\mathbf{X}^n=\mathbf{A}\end{equation} with $n\geq 2, n\;\epsilon\;\mathbb{N} $.
The way to the solution consists of two steps:\\
1. Transform the quaternion equation into the  isomorphic matrix equation.
Instead of the equation $\mathbf{X}^n=\mathbf{A}$ mit $\mathbf{A}=\quat{a}=\vekk{a}$ we solve the equation
 \begin{equation}\label{xhochn}\mathfrak{X}^n=\mathfrak{A}\end{equation}
 with $\mathfrak{A}=\qmatrix{a}$.\\
 The solution set is all matrices $\mathfrak{X}$, satisfying this equation.\\
 2. From the solution set $\{\mathfrak{X}\}$ we select the subset\\
  $\mathfrak{X}_Q=\{\mathfrak{X}\|\;
  \mathfrak{X}\;\mbox{isomorphic to a quaternion} \}$.\\
The isomorphic relation provides the solution set of the quaternion equation directly.\\
If we want to solve equation (\ref{xhochn}), then we only have to transcribe the solution from a good textbook
about the theory of matrices. Therefore we regard at first the structur of the matrices which are isomorphic
to quaternion (we call they q-matrices in the following). This preliminary work is necessary, because we need
the Jordan normal form of the q-matrix $\mathfrak{A}$.
\section{Jordan normal form of q-matrices}
All of the occuring numbers $a_i, x_i$ etc. are komplex numbers (elements of $\Bbb{C}$).\\
Solving problem: Compute the transformation matrix $\mathfrak{U}$ and the Jordan matrix $\mathfrak{J}$ in the
equation $\mathfrak{A}=\mathfrak{U}\mathfrak{J}\mathfrak{U}^{-1}$ for a given q-matrix $\mathfrak{A}$.
Q-matrices are of the form:
\[\mathfrak{A}=\qmatrix{a}\]
The characteristic polynomial of the q-matrix $\mathfrak{A}$ is $c(y)=(y^2-2a_0y+a_0^2+\mathbf{a}^2)^2$.\\
The minimal polynomial of $\mathfrak{A}$ is $m(y)=y^2-2a_0y+a_0^2+\mathbf{a}^2$ and the reduced adjugate matrix
\begin{equation}\label{redadj}\mathfrak{C}(y)=\admatrixy{a}{y}\end{equation}
The solutions of the characteristic polynomial are
\begin{equation}\label{nullst}y_1=y_2=a_0+i\sqrt{\mathbf{a}^2},\quad y_3=y_4=a_0-i\sqrt{\mathbf{a}^2}\end{equation}
There are three cases for the Jordan normal form of $\mathfrak{A}=\mathfrak{U}\mathfrak{J}\mathfrak{U}^{-1}$:\\
1. $\mathbf{a}=\mathbf{0}$.\\
In this case $\mathfrak{A}$ is of the representation
\begin{equation}\label{J4}\mathfrak{A}=\mathfrak{J}=
\left(\begin{array}{*{3}{c@{\;}}c}
	        a_0 & 0 & 0 & 0\\
			0 &  a_0 & 0 &  0\\
			0 &  0 &  a_0 & 0\\
			0 & 0 &  0 &  a_0
      \end{array}\right)
	  = a_0\left(\begin{array}{*{3}{c@{\;}}c}
	        1 & 0 & 0 & 0\\
			0 &  1 & 0 &  0\\
			0 &  0 &  1 & 0\\
			0 & 0 &  0 &  1
      \end{array}\right)\end{equation}	
The matrix $\mathfrak{A}$ already has Jordan normal form with (after reducing) the minimal polynomial
$m(y)= y-a_0$ and the reduced adjugate matrix
\[\mathfrak{C}(y)=\mathfrak{E}=\left(\begin{array}{*{3}{c@{\;}}c}
	        1 & 0 & 0 & 0\\
			0 &  1 & 0 &  0\\
			0 &  0 &  1 & 0\\
			0 & 0 &  0 &  1
      \end{array}\right)\]
2. $\mathbf{a}\neq\mathbf{0}$ und $\mathbf{a}^2\neq 0$\\
Here the minimal polynomial has 2 simple roots. The corresponding Jordan matrix is diagonalizable:
\begin{equation}\label{Jd}\mathfrak{J}=\left(\begin{array}{*{3}{c@{\;}}c}
	        a_0+i\sqrt{\mathbf{a}^2} & 0 & 0 & 0\\
			0 &  a_0+i\sqrt{\mathbf{a}^2} & 0 &  0\\
			0 &  0 &  a_0-i\sqrt{\mathbf{a}^2} & 0\\
			0 & 0 &  0 &  a_0-i\sqrt{\mathbf{a}^2}
      \end{array}\right)\end{equation}
The transformation matrix $\mathfrak{U}$ is assembled by the columns of the reduced adjugate matrix
 $\mathfrak{C}(y)$ (\ref{redadj}).
Because $\mathbf{a}\neq\mathbf{0}$ and $\mathbf{a}^2\neq 0$, we have to differentiate between three cases:\\
 a) $a_2^2+a_3^2\neq 0$\\
$\mathfrak{U}$ is assembled by the first both columns of
$\mathfrak{C}(a_0+i\sqrt{\mathbf{a}^2})$ and the last both columns of $\mathfrak{C}(a_0-i\sqrt{\mathbf{a}^2})$
Then we have:
\begin{equation}\label{trans11}\mathfrak{U}=
\left(\begin{array}{*{3}{c@{\;}}c}
i\sqrt{\mathbf{a}^2} & -a_1 & -a_2 & -a_3\\
a_1 & i\sqrt{\mathbf{a}^2} & -a_3 & a_2\\
a_2 & a_3 & -i\sqrt{\mathbf{a}^2} & -a_1\\
a_3 & -a_2 & a_1 & -i\sqrt{\mathbf{a}^2}
\end{array}\right)\end{equation} and
\[\mathfrak{U}^{-1}=\frac{1}{2(a_2^2+a_3^2)\sqrt{\mathbf{a}^2}}\;\cdot\]
\[\left(\begin{array}{*{3}{c@{\:}}c}
-i(a_2^2+a_3^2) & 0 & a_2\sqrt{\mathbf{a}^2}-ia_1a_3 & a_3\sqrt{\mathbf{a}^2}+ia_1a_2\\
0 & -i(a_2^2+a_3^2) & a_3\sqrt{\mathbf{a}^2}+ia_1a_2 & -a_2\sqrt{\mathbf{a}^2}+ia_1a_3\\
-a_2\sqrt{\mathbf{a}^2}-ia_1a_3 & -a_3\sqrt{\mathbf{a}^2}+ia_1a_2 & i(a_2^2+a_3^2) & 0\\
-a_3\sqrt{\mathbf{a}^2}+ia_1a_2 & a_2\sqrt{\mathbf{a}^2}+ia_1a_3 & 0 & i(a_2^2+a_3^2)
\end{array}\right)\]
b) $a_1^2+a_3^2\neq 0$\\
$\mathfrak{U}$ results from the first and third column of
$\mathfrak{C}(a_0+i\sqrt{\mathbf{a}^2})$ and the second and fourth column of
$\mathfrak{C}(a_0-i\sqrt{\mathbf{a}^2})$
Then we have:
\begin{equation}\label{trans12}\mathfrak{U}=
\left(\begin{array}{*{3}{c@{\;}}c}
i\sqrt{\mathbf{a}^2} & -a_2 & -a_1 & -a_3\\
a_1 & -a_3 & -i\sqrt{\mathbf{a}^2} & a_2\\
a_2 & i\sqrt{\mathbf{a}^2} & a_3 & -a_1\\
a_3 & a_1 & -a_2 & -i\sqrt{\mathbf{a}^2}
\end{array}\right)\end{equation} and\\
\[\mathfrak{U}^{-1}=\frac{1}{2(a_1^2+a_3^2)\sqrt{\mathbf{a}^2}}\;\cdot\]
\[\left(\begin{array}{*{3}{c@{\:}}c}
-i(a_1^2+a_3^2) & a_1\sqrt{\mathbf{a}^2}+ia_2a_3 & 0 & a_3\sqrt{\mathbf{a}^2}-ia_1a_2\\
0 & -a_3\sqrt{\mathbf{a}^2}+ia_1a_2 & -i(a_1^2+a_3^2) & a_1\sqrt{\mathbf{a}^2}+ia_2a_3\\
-a_1\sqrt{\mathbf{a}^2}+ia_2a_3 & +i(a_1^2+a_3^2) & a_3\sqrt{\mathbf{a}^2}+ia_1a_2 & 0\\
-a_3\sqrt{\mathbf{a}^2}-ia_1a_2 & 0 & -a_1\sqrt{\mathbf{a}^2}+ia_2a_3 & i(a_1^2+a_3^2)
\end{array}\right)\]
c) $a_1^2+a_2^2\neq 0$\\
$\mathfrak{U}$ results from the first and fourth column of
$\mathfrak{C}(a_0+i\sqrt{\mathbf{a}^2})$ and the second and third column of
$\mathfrak{C}(a_0-i\sqrt{\mathbf{a}^2})$
Then we have:
\begin{equation}\label{trans13}\mathfrak{U}=
\left(\begin{array}{*{3}{c@{\;}}c}
i\sqrt{\mathbf{a}^2} & -a_3 & -a_1 & -a_2\\
a_1 & a_2 & -i\sqrt{\mathbf{a}^2} & -a_3\\
a_2 & -a_1 & a_3 & -i\sqrt{\mathbf{a}^2}\\
a_3 & i\sqrt{\mathbf{a}^2} & -a_2 & a_1
\end{array}\right)\end{equation} and
\[\mathfrak{U}^{-1}=\frac{1}{2(a_1^2+a_2^2)\sqrt{\mathbf{a}^2}}\;\cdot\]
\[\left(\begin{array}{*{3}{c@{\:}}c}
-i(a_1^2+a_2^2) & a_1\sqrt{\mathbf{a}^2}-ia_2a_3 & a_2\sqrt{\mathbf{a}^2}+ia_1a_3 & 0\\
0 & a_2\sqrt{\mathbf{a}^2}+ia_1a_3 & -a_1\sqrt{\mathbf{a}^2}+ia_2a_3 & -i(a_1^2+a_2^2)\\
-a_1\sqrt{\mathbf{a}^2}-ia_2a_3 & i(a_1^2+a_2^2) & 0 & -a_2\sqrt{\mathbf{a}^2}+ia_1a_3\\
-a_2\sqrt{\mathbf{a}^2}+ia_1a_3 & 0 & i(a_1^2+a_3^2) & a_1\sqrt{\mathbf{a}^2}+ia_2a_3
\end{array}\right)\]
3. $\mathbf{a}\neq\mathbf{0}$ und $\mathbf{a}^2 = 0$\\
Here the minimal polynomial is  $m(y)=y^2-2a_0y+a_0^2=(y-a_0)^2$.
The only (of multiplicity 2) solution is $y_1=a_0$.
The corresponding Jordan matrix of $\mathfrak{A}$ is then:
\begin{equation}\label{Jk}\mathfrak{J}=\left(\begin{array}{*{3}{c@{\;}}c}
	        a_0 & 1 & 0 & 0\\
			0 &  a_0 & 0 &  0\\
			0 &  0 &  a_0 & 1\\
			0 & 0 &  0 &  a_0
           \end{array}\right)\end{equation}
We use the reduced adjugate matrix $\mathfrak{C}(y)$ and the derivative:
\[\mathfrak{C}(y)=\admatrixy{a}{y}\]
\[\mathfrak{C}\;'(y)=\Ematrix\]
to compute the transformation matrix $\mathfrak{U}$.
We have to differentiate between three cases too:\\
a) $a_2^2+a_3^2\neq 0$\\
$\mathfrak{U}$ results from the first and the second column of  $\mathfrak{C}(a_0)$ and of $\mathfrak{C}\;'(a_0)$.
\begin{equation}\label{trans21}\mathfrak{U}=
\left(\begin{array}{*{3}{c@{\;\;}}c}
0 & 1 & -a_1 & 0\\
a_1 & 0 & 0 & 1\\
a_2 & 0 & a_3 & 0\\
a_3 & 0 & -a_2 & 0
\end{array}\right)\end{equation} and
\[\mathfrak{U}^{-1}=\frac{1}{a_2^2+a_3^2}\;\cdot
\left(\begin{array}{*{3}{c@{\;\;}}c}
0 & 0 & a_2 & a_3\\
a_2^2+a_3^2 & 0 & a_1a_3 & -a_1a_2\\
0 & 0 & a_3 & -a_2\\
0 & a_2^2+a_3^2 & -a_1a_2 & -a_1a_3
\end{array}\right)\]
b) $a_1^2+a_3^2\neq 0$\\
$\mathfrak{U}$ results from the first and the third column of $\mathfrak{C}(a_0)$ and of $\mathfrak{C}\;'(a_0)$.
\begin{equation}\label{trans22}\mathfrak{U}=
\left(\begin{array}{*{3}{c@{\;\;}}c}
0 & 1 & -a_2 & 0\\
a_1 & 0 & -a_3 & 0\\
a_2 & 0 & 0 & 1\\
a_3 & 0 & a_1 & 0
\end{array}\right)\end{equation} and
\[\mathfrak{U}^{-1}=\frac{1}{a_1^1+a_3^2}\;\cdot
\left(\begin{array}{*{3}{c@{\:\:}}c}
0 & a_1 & 0 & a_3\\
a_1^2+a_3^2 & -a_2a_3 & 0 & a_1a_2\\
0 & -a_3 & 0 & a_1\\
0 & -a_1a_2 & a_1^2+a_3^2 & -a_2a_3
\end{array}\right)\]
c) $a_1^2+a_2^2\neq 0$\\
$\mathfrak{U}$ results from the first and the fourth column of $\mathfrak{C}(a_0)$ and of $\mathfrak{C}\;'(a_0)$.
\begin{equation}\label{trans23}\mathfrak{U}=
\left(\begin{array}{*{3}{c@{\;}}c}
0 & 1 & -a_3 & 0\\
a_1 & 0 & a_2 & 0\\
a_2 & 0 & -a_1 & 0\\
a_3 & 0 & 0 & 1
\end{array}\right)\end{equation} and
\[\mathfrak{U}^{-1}=\frac{1}{a_1^1+a_2^2}\;\cdot
\left(\begin{array}{*{3}{c@{\:\:}}c}
0 & a_1 & a_2 & 0\\
a_1^2+a_2^2 & a_2a_3 & -a_1a_3 & 0\\
0 & a_2 & -a_1 & 0\\
0 & -a_1a_3 & -a_2a_3 & a_1^2+a_2^2
\end{array}\right)\]
Why do I compute the transformation matrix $\mathfrak{U}$ and the Jordan matrix $\mathfrak{J}$ so extensive?
In order to show that $\mathfrak{U}$ and $\mathfrak{J}$ are uniquely determined by the four components of the
quaternion $\mathbf{A}$. Futhermore we found that a q-matrix is similar to one of three Jordan matrices
(\ref{J4}), (\ref{Jd}) or (\ref{Jk}). Now we can go about the solution of (\ref{xhochn}).
For it we remove the appropriate chapter from a suitable textbook.
\section{Solution}
We also consult  \cite{ga} chapter 8.6 at page 243. There we read that we have to make yet another case
differentiation, the differentiation between invertible and singular $\mathfrak{A}$.
Now let us begin with
\subsection{Invertible matrices/quaternions}
1. Case - The matrix $\mathfrak{A}$ is of type (\ref{J4}).\\
It can be written as
\begin{equation}\label{Formel1}\mathfrak{A}= a_0\left(\begin{array}{*{3}{c@{\;}}c}
	        1 & 0 & 0 & 0\\
			0 &  1 & 0 &  0\\
			0 &  0 &  1 & 0\\
			0 & 0 &  0 &  1
      \end{array}\right)\end{equation}
Furthermore we need all invertible matrices $\mathfrak{V}$ (chapter 8.6 from \cite{ga}) which are commutable with
$\mathfrak{A}$ also satisfy the equation $\mathfrak{A}\mathfrak{V}=\mathfrak{V}\mathfrak{A}$.
Because $\mathfrak{A}=a_0\mathfrak{E}$ the above equation holds for any invertible matrix $\mathfrak{V}$.
Let $w_i$ be an arbitrary n-th root of $a_0$. The general solution of
(\ref{xhochn}) is then:
\begin{equation}\label{Formel1l}\mathfrak{X}= \mathfrak{V}\left(\begin{array}{*{3}{c@{\;}}c}
	        w_1 & 0 & 0 & 0\\
			0 &  w_2 & 0 &  0\\
			0 &  0 &  w_3 & 0\\
			0 & 0 &  0 &  w_4
      \end{array}\right)\mathfrak{V}^{-1}\end{equation}
Here only the solution matrices are of interest, which are isomorphic to a quaternion. Than there are the the
following cases:\\
a) $w_1=w_2=w_3=w_4$:  Then is $\mathfrak{X}=\mathfrak{V}w_1\mathfrak{E}\mathfrak{V}^{-1}=w_1\mathfrak{E}$.\\
Solution quaternions $\mathbf{X}$ are any $\mathbf{X}=w_1\einsq$. There are $n$ solutions according to the
$n$ different values of $\sqrt[n]{a_0}$ (This is the trivial solution.).\\[1ex]
b) $w_1=w_2$ und $w_3=w_4$ \\
Here $w_1$ and $w_3$ are two different n-th roots of $a_0$, t. i.
$w_1^n=w_3^n=a_0$ but $w_1\neq w_3$. So that $\mathfrak{X}$ is isomorphic to a quaternion, we need to construct
$\mathfrak{U}$ according to  (\ref{trans11}), (\ref{trans12}) or (\ref{trans13}). If we select any numbers for
$x_0-x_4$ then we can compute the $\mathfrak{U}$ without restrictions.
The general solution is
\begin{equation}\mathfrak{X}= \mathfrak{V}\left(\begin{array}{*{3}{c@{\;}}c}
	        w_1 & 0 & 0 & 0\\
			0 &  w_1 & 0 &  0\\
			0 &  0 &  w_3 & 0\\
			0 & 0 &  0 &  w_3
      \end{array}\right)\mathfrak{V}^{-1}\end{equation}
with arbitrary invertible matrix $\mathfrak{V}$.
Therewith $\mathfrak{X}$ is isomorphic to a quaternion (of type (\ref{Jd})) we have the conditions
\begin{equation}\label{Xwert1}\begin{array}{c}
         w_1=x_0+i\sqrt{x^2}\\
		 w_3=x_0-i\sqrt{x^2}
\end{array}\end{equation}
and the matrix $\mathfrak{V}$ has to be an U-matrix according to (\ref{trans11}), (\ref{trans12}) or
(\ref{trans13}).
Because $w_1$ and $w_3$ are $\neq 0$ we always can compute the corresponding u-matrix and $x_0 ... x_4$ only need
to satisfy  (\ref{Xwert1}).
The corresponding quaternion set $\{\mathbf{X}\}$ results into
\begin{equation}
\{\mathbf{X}\}=\{\mathbf{X}\mid x_0=\frac{1}{2}(w_1+w_3) \wedge
\mathbf{x}^2=x_1^2+x_2^2+x_3^2=-\frac{1}{4}(w_1-w_3)^2\}
\end{equation}
There are ${n \choose 2}$ of this solution sets (accordingly to the numbers of the possible combinations of $w_1$
and $w_3$).\\

1. Example: n=2, to solve the equation $\mathbf{X}^2=a_0\mathbf{E}=\ixsqv{a_0}$.
The simple solution is
$\mathbf{X}=\pm\sqrt{a_0}\einsqv$. The second solution (according to b))is of more interest. The two different
roots $w_1$ and $w_3$ are: $w_1=\sqrt{a_0}$ and $w_3=-\sqrt{a_0}$. So that the solution set is
\{$\mathbf{X}$\} with $x_0=0$ and $\mathbf{x}^2=-a_0$\\
2. Example: n=3, to solve the equation $\mathbf{X}^3=a_0\mathbf{E}=\ixsqv{a_0}$. The trivial solution is
obviously $\mathbf{X}=\sqrt[3]{a_0}\;\mathbf{E}$, whereas for $\sqrt[3]{a_0}$ tree values are possible.
The solutions according to b) are ${3 \choose 2}=3$ different sets, which are computed in the following. The 3
different 3-rd roots from $a_0$ are: $w_1=\sqrt[3]{a_0}$,
$w_2=w_1(-\frac 1 2 + \frac i 2 \sqrt{3})$ and $w_3=w_1(-\frac 1 2 - \frac i 2 \sqrt{3})$.
Then the  combinations result: ($w_1, w_2$), ($w_1, w_3$) and ($w_2, w_3$). This provides the following three
solution sets:\\
($w_1, w_2$)$\;\Rightarrow\; x_0=w_1(\frac 1 4 +\frac i 4 \sqrt{3})\quad \mathbf{x}^2=\frac 3 4 w_1^2
(-\frac 1 2 + \frac i 2 \sqrt{3})$\\
($w_1, w_3$)$\;\Rightarrow\;x_0=w_1(\frac 1 4 -\frac i 4 \sqrt{3})\quad \mathbf{x}^2=\frac 3 4 w_1^2
(-\frac 1 2 - \frac i 2 \sqrt{3})$\\
($w_2, w_3$)$\;\Rightarrow\;x_0=-\frac 1 2 w_1\quad \mathbf{x}^2=\frac 3 4 w_1^2$. \\[1.5ex]
2. Case: The matrix $\mathfrak{A}$ is of type (\ref{Jd}).\\
Using a transformation matrix of kind (\ref{trans11}), (\ref{trans12}) or (\ref{trans13}) the $\mathfrak{A}$
can be written as: \quad$\mathfrak{A}=\mathfrak{U}\mathfrak{J}\mathfrak{U}^{-1}$.\\
$\mathfrak{J}$ is a diagonal matrix. An arbitrary, with $\mathfrak{J}$ commutable, invertible matrix
$\mathfrak{V}$ is of kind: (q.v. \cite{ga}, chapter 8.2)
\begin{equation}\label{Matrix5}\mathfrak{V}= \left(\begin{array}{*{3}{c@{\;}}c}
	v_1 & v_2 & 0 & 0 \\
	v_3 & v_4 & 0 & 0\\
	0 & 0 & v_5 & v_6\\
	0 & 0 & v_7 & v_8
\end{array}\right)\end{equation}
Then the general solution of (\ref{xhochn}) is (q.v. \cite{ga}, S. 245)
\begin{equation}\mathfrak{X}= \mathfrak{U}\mathfrak{V}\left(\begin{array}{*{3}{c@{\;}}c}
	        \sqrt[n]{a_0+i\sqrt{\mathbf{a}^2}} & 0 & 0 & 0\\
			0 &  \sqrt[n]{a_0+i\sqrt{\mathbf{a}^2}} & 0 &  0\\
			0 &  0 &  \sqrt[n]{a_0-i\sqrt{\mathbf{a}^2}} & 0\\
			0 & 0 &  0 &  \sqrt[n]{a_0-i\sqrt{\mathbf{a}^2}}
      \end{array}\right)\mathfrak{V}^{-1}\mathfrak{U}^{-1}\end{equation}
In this general formula different values can be inserted for the n-th root in the 1. and 2. as soon as in the
3. and 4. row (accordingly the multiplicity  of the n-th root). We need only $\mathfrak{X}$-matrices which are
isomorphic to a quaternion. Therefore roots of the 1. and 2. row as soon as the 3. and 4. row of $\mathfrak{J}$
have to be equal respectively. But then
 \begin{equation}\mathfrak{J}_x=\left(\begin{array}{*{3}{c@{\;}}c}
 		\sqrt[n]{a_0+i\sqrt{\mathbf{a}^2}} & 0 & 0 & 0\\
			0 &  \sqrt[n]{a_0+i\sqrt{\mathbf{a}^2}} & 0 &  0\\
			0 &  0 &  \sqrt[n]{a_0-i\sqrt{\mathbf{a}^2}} & 0\\
			0 & 0 &  0 &  \sqrt[n]{a_0-i\sqrt{\mathbf{a}^2}}
      \end{array}\right)\end{equation}
is commutable with $\mathfrak{V}$ and the equation
$\mathfrak{J}_x=\mathfrak{V}\mathfrak{J}_x\mathfrak{V}^{-1}$ is true.
$\mathfrak{X}$ is:
\begin{equation}\mathfrak{X}= \mathfrak{U}\left(\begin{array}{*{3}{c@{\;}}c}
	        \sqrt[n]{a_0+i\sqrt{\mathbf{a}^2}} & 0 & 0 & 0\\
			0 &  \sqrt[n]{a_0+i\sqrt{\mathbf{a}^2}} & 0 &  0\\
			0 &  0 &  \sqrt[n]{a_0-i\sqrt{\mathbf{a}^2}} & 0\\
			0 & 0 &  0 &  \sqrt[n]{a_0-i\sqrt{\mathbf{a}^2}}
      \end{array}\right)\mathfrak{U}^{-1}.\end{equation}
The transformation matrix $\mathfrak{U}$ is as same as that of $\mathfrak{A}$. Then we can write $\mathbf{X}$ as
quaternion as follows:
\begin{equation}\label{Lösung2}
\mathbf{X}=\sqrt[n]{a_0+i\sqrt{\mathbf{a}^2}}\left(
	\begin{array}{c}\frac{1}{2}\\[1ex]
	\frac{1}{2i}\frac{\mathbf{a}}{\sqrt{\mathbf{a}^2}}
	\end{array}\right)+\sqrt[n]{a_0-i\sqrt{\mathbf{a}^2}}\left(
	\begin{array}{c}\frac{1}{2}\\[1ex]
	-\frac{1}{2i}\frac{\mathbf{a}}{\sqrt{\mathbf{a}^2}}
	\end{array}\right).\end{equation}
n values can be inserted for the n-th roots so that follows alltogether $n^2$ solutions.	
(every value of $\sqrt[n]{a_0+i\sqrt{\mathbf{a}^2}}$ can be combined with every value of
$\sqrt[n]{a_0-i\sqrt{\mathbf{a}^2}}$.) here an excample:\\
Let $\mathbf{A}=\left(\begin{array}{c}76\\ 24i\\ 12i\\ -41i\end{array}\right)$, $n=3$, to solve the equation
$\mathbf{X}^3=\mathbf{A}$. The two eigenvalues of $\mathfrak{A}$ are with $\sqrt{\mathbf{a}^2}=49i$:\\
 $a_0+i\sqrt{\mathbf{a}^2}=27$ and $a_0-i\sqrt{\mathbf{a}^2}=125$. The cubic roots have the values:
 $3;\; -\frac{3}{2}+\frac{3}{2}i\sqrt{3};\;-\frac{3}{2}-\frac{3}{2}i\sqrt{3}$ as soon as
 $5;\; -\frac{5}{2}+\frac{5}{2}i\sqrt{3};\;-\frac{5}{2}-\frac{5}{2}i\sqrt{3}$. The 9 solutions are written down
 as:\\
$\mathbf{X}_1=3\left(\begin{array}{c}\frac{1}{2}\\[.5ex]\frac{-12i}{49}\\[.5ex]\frac{-6i}{49}\\[.5ex]
\frac{41i}{98}\end{array}\right)+5\left(\begin{array}{c}\frac{1}{2}\\[.5ex]\frac{12i}{49}\\[.5ex]\frac{6i}{49}\\
[.5ex]\frac{-41i}{98}\end{array}\right)=\left(\begin{array}{c}4\\[.5ex]\frac{24i}{49}\\[.5ex]\frac{12i}{49}\\[.5ex]
\frac{-41i}{49}\end{array}\right)$\\[1ex]
$\mathbf{X}_2=(-\frac{3}{2}+\frac{3}{2}i\sqrt{3})\left(\begin{array}{c}\frac{1}{2}\\[.5ex]\frac{-12i}{49}\\[.5ex]
\frac{-6i}{49}\\[.5ex]\frac{41i}{98}\end{array}\right)+5\left(\begin{array}{c}\frac{1}{2}\\[.5ex]\frac{12i}{49}\\
[.5ex]\frac{6i}{49}\\[.5ex]\frac{-41i}{98}\end{array}\right)=\left(\begin{array}{c}\frac{7+3i\sqrt{3}}{4}\\[.5ex]
\frac{18\sqrt{3}+78i}{49}\\[.5ex]\frac{9\sqrt{3}+39i}{49}\\[.5ex]\frac{-123\sqrt{3}-533i}{196}
\end{array}\right)$\\[1ex]
$\mathbf{X}_3=(-\frac{3}{2}-\frac{3}{2}i\sqrt{3})\left(\begin{array}{c}\frac{1}{2}\\[.5ex]\frac{-12i}{49}\\[.5ex]
\frac{-6i}{49}\\[.5ex]\frac{41i}{98}\end{array}\right)+5\left(\begin{array}{c}\frac{1}{2}\\[.5ex]\frac{12i}{49}\\
[.5ex]\frac{6i}{49}\\[.5ex]\frac{-41i}{98}\end{array}\right)=\left(\begin{array}{c}\frac{7-3i\sqrt{3}}{4}\\[.5ex]
\frac{-18\sqrt{3}+78i}{49}\\[.5ex]\frac{-9\sqrt{3}+39i}{49}\\[.5ex]\frac{123\sqrt{3}-533i}{196}
\end{array}\right)$\\[1ex]
$\mathbf{X}_4=3\left(\begin{array}{c}\frac{1}{2}\\[.5ex]\frac{-12i}{49}\\[.5ex]\frac{-6i}{49}\\[.5ex]
\frac{41i}{98}\end{array}\right)+(-\frac{5}{2}+\frac{5}{2}i\sqrt{3})\left(\begin{array}{c}\frac{1}{2}\\[.5ex]
\frac{12i}{49}\\[.5ex]\frac{6i}{49}\\[.5ex]\frac{-41i}{98}\end{array}\right)=\left(\begin{array}{c}
\frac{1+5i\sqrt{3}}{4}\\[.5ex]\frac{-30\sqrt{3}-78i}{49}\\[.5ex]\frac{-15\sqrt{3}-33i}{49}\\[.5ex]
\frac{205\sqrt{3}+451i}{196}\end{array}\right)$\\[1ex]
$\mathbf{X}_5=(-\frac{3}{2}+\frac{3}{2}i\sqrt{3})\left(\begin{array}{c}\frac{1}{2}\\[.5ex]\frac{-12i}{49}\\[.5ex]
\frac{-6i}{49}\\[.5ex]\frac{41i}{98}\end{array}\right)+(-\frac{5}{2}+\frac{5}{2}i\sqrt{3})\left(\begin{array}{c}
\frac{1}{2}\\[.5ex]\frac{12i}{49}\\[.5ex]\frac{6i}{49}\\[.5ex]\frac{-41i}{98}\end{array}\right)=\left(
\begin{array}{c}-2+2i\sqrt{3}\\[.5ex]\frac{-12\sqrt{3}-12i}{49}\\[.5ex]\frac{-6\sqrt{3}-6i}{49}\\[.5ex]
\frac{41\sqrt{3}+41i}{98}\end{array}\right)$\\[1ex]
$\mathbf{X}_6=(-\frac{3}{2}-\frac{3}{2}i\sqrt{3})\left(\begin{array}{c}\frac{1}{2}\\[.5ex]\frac{-12i}{49}\\[.5ex]
\frac{-6i}{49}\\[.5ex]\frac{41i}{98}\end{array}\right)+(-\frac{5}{2}+\frac{5}{2}i\sqrt{3})\left(\begin{array}{c}
\frac{1}{2}\\[.5ex]\frac{12i}{49}\\[.5ex]\frac{6i}{49}\\[.5ex]\frac{-41i}{98}\end{array}\right)=\left(
\begin{array}{c}-2+\frac{1}{2}i\sqrt{3}\\[.5ex]\frac{-48\sqrt{3}-12i}{49}\\[.5ex]\frac{-24\sqrt{3}-6i}{49}\\[.5ex]
\frac{164\sqrt{3}+41i}{98}\end{array}\right)$\\[1ex]
$\mathbf{X}_7=3\left(\begin{array}{c}\frac{1}{2}\\[.5ex]\frac{-12i}{49}\\[.5ex]\frac{-6i}{49}\\[.5ex]
\frac{41i}{98}\end{array}\right)+(-\frac{5}{2}-\frac{5}{2}i\sqrt{3})\left(\begin{array}{c}\frac{1}{2}\\[.5ex]
\frac{12i}{49}\\[.5ex]\frac{6i}{49}\\[.5ex]\frac{-41i}{98}\end{array}\right)=\left(\begin{array}{c}
\frac{1-5i\sqrt{3}}{4}\\[.5ex]\frac{30\sqrt{3}-78i}{49}\\[.5ex]\frac{15\sqrt{3}-33i}{49}\\[.5ex]
\frac{-205\sqrt{3}+451i}{196}\end{array}\right)$\\[1ex]
$\mathbf{X}_8=(-\frac{3}{2}+\frac{3}{2}i\sqrt{3})\left(\begin{array}{c}\frac{1}{2}\\[.5ex]\frac{-12i}{49}\\[.5ex]
\frac{-6i}{49}\\[.5ex]\frac{41i}{98}\end{array}\right)+(-\frac{5}{2}-\frac{5}{2}i\sqrt{3})\left(\begin{array}{c}
\frac{1}{2}\\[.5ex]\frac{12i}{49}\\[.5ex]\frac{6i}{49}\\[.5ex]\frac{-41i}{98}\end{array}\right)=\left(
\begin{array}{c}-2-\frac{1}{2}i\sqrt{3}\\[.5ex]\frac{48\sqrt{3}-12i}{49}\\[.5ex]\frac{24\sqrt{3}-6i}{49}\\[.5ex]
\frac{-164\sqrt{3}+41i}{98}\end{array}\right)$\\[1ex]
$\mathbf{X}_9=(-\frac{3}{2}-\frac{3}{2}i\sqrt{3})\left(\begin{array}{c}\frac{1}{2}\\[.5ex]\frac{-12i}{49}\\[.5ex]
\frac{-6i}{49}\\[.5ex]\frac{41i}{98}\end{array}\right)+(-\frac{5}{2}-\frac{5}{2}i\sqrt{3})\left(\begin{array}{c}
\frac{1}{2}\\[.5ex]\frac{12i}{49}\\[.5ex]\frac{6i}{49}\\[.5ex]\frac{-41i}{98}\end{array}\right)=\left(
\begin{array}{c}-2-2i\sqrt{3}\\[.5ex]\frac{12\sqrt{3}-12i}{49}\\[.5ex]\frac{6\sqrt{3}-6i}{49}\\[.5ex]
\frac{-41\sqrt{3}+41i}{98}\end{array}\right)$\\[1ex]
3. Case: The matrix $\mathfrak{A}$ is of type (\ref{Jk}).\\
Here the transformation matrix is of kind (\ref{trans21}), (\ref{trans22}) or (\ref{trans23}) and
$\mathfrak{J}$ is of kind (\ref{Jk}). We have the equation
$\mathfrak{A}=\mathfrak{U}\mathfrak{J}\mathfrak{U}^{-1}$. An arbitrary, invertible, with $\mathfrak{J}$
commutable matrix $\mathfrak{V}$ has the form (according to \cite{ga}, Kapitel 8.2):
\begin{equation}\label{Matrix6}\mathfrak{V}= \left(\begin{array}{*{3}{c@{\;}}c}
	v_1 & v_2 & v_3 & v_4 \\
	0 & v_1 & 0 & v_3\\
	v_5 & v_6 & v_7 & v_8\\
	0 & v_5 & 0 & v_7
\end{array}\right)\end{equation}
The general solution of (\ref{xhochn}) is then (with the abbreviation $w=\sqrt[n]{a_0}$)
\begin{equation}\mathfrak{X}= \mathfrak{U}\mathfrak{V}\left(\begin{array}{*{3}{c@{\;}}c}
	        w_1 & \frac{1}{nw_1^{n-1}} & 0 & 0\\
			0 &  w_1 & 0 &  0\\
			0 &  0 &  w_2 & \frac{1}{nw_2^{n-1}}\\
			0 & 0 &  0 &  w_2
      \end{array}\right)\mathfrak{V}^{-1}\mathfrak{U}^{-1}\end{equation}
We need only the solutions $\mathfrak{X}$, which are isomorphic to a quaternion. In this case the two roots
$w_1$ and $w_2$ have to be equal. It follows the solution matrix:
\begin{equation}\mathfrak{X}= \mathfrak{U}\mathfrak{V}\left(\begin{array}{*{3}{c@{\;}}c}
	        w & \frac{1}{nw^{n-1}} & 0 & 0\\
			0 &  w & 0 &  0\\
			0 &  0 &  w & \frac{1}{nw^{n-1}}\\
			0 & 0 &  0 &  w
      \end{array}\right)\mathfrak{V}^{-1}\mathfrak{U}^{-1}\end{equation}
In this formula the w-matrix  is commutable with $\mathfrak{V}$ too and $\mathfrak{X}$ reduces to
\begin{equation}\mathfrak{X}= \mathfrak{U}\left(\begin{array}{*{3}{c@{\;}}c}
	        w & \frac{1}{nw^{n-1}} & 0 & 0\\
			0 &  w & 0 &  0\\
			0 &  0 &  w & \frac{1}{nw^{n-1}}\\
			0 & 0 &  0 &  w
      \end{array}\right)\mathfrak{U}^{-1}\end{equation}
The retranslation into the quaternion notation provides then
\[\mathbf{X}=w\einsq +\frac{1}{nw^{n-1}}\quatvek{a}=w\einsq +\frac{w}{na_0}\quatvek{a}\]
\begin{equation}\label{Lösung3}\mathbf{X}=w\left(\begin{array}{c}1\\[1ex]
\displaystyle{\frac{1}{na_0}\mathbf{a}}\end{array}\right)
\end{equation}
Thereby we have $w=\sqrt[n]{a_0}$ and there are $n$ different values possible.\\
Example: Let $\mathbf{A}=\left(\begin{array}{c}1\\ 3i\\ 4i\\ 5\end{array}\right)$, $n=4$.
Then it follows $w=\sqrt[4]{1}$ with the solution set $\{1;-1;i;-i\}$. The 4 solutions of $\mathbf{X}^4=\mathbf{A}$
are:
\[\mathbf{X}_1=\left(\begin{array}{c}1\\ \frac{3}{4}i\\ i\\ \frac{5}{4}\end{array}\right)\quad
\mathbf{X}_2=\left(\begin{array}{c}-1\\ -\frac{3}{4}i\\ -i\\ -\frac{5}{4}\end{array}\right)\]
\[\mathbf{X}_3=\left(\begin{array}{c}i\\ -\frac{3}{4}\\ -1\\ \frac{5}{4}i\end{array}\right)\quad
\mathbf{X}_4=\left(\begin{array}{c}-i\\ \frac{3}{4}\\ 1\\ -\frac{5}{4}i\end{array}\right)\]
\subsection{Singular matrices/quaternions}
Singular quaternions can too be divided into the three classes (\ref{J4}), (\ref{Jd}) and (\ref{Jk}).
For $\mathfrak{A}$ we can make the same case differentiation as above at the invertible quaternions.
Then we get different solutions as in the case of invertible quaternions.\\
1. $\mathbf{a}=\mathbf{0}$. (Type (\ref{J4}))\\
So that $\mathbf{A}$, resp. $\mathfrak{A}$ is the Zero quaternion, resp. the Zero matrix.
The following equation is to solve:
\[\mathfrak{X}^n=\left(\begin{array}{*{3}{c@{\;}}c}
	0 & 0 & 0 & 0 \\
	0 & 0 & 0 & 0\\
	0 & 0 & 0 & 0\\
	0 & 0 & 0 & 0
\end{array}\right)\]
Evidently the Zero matrix is commutable with any matrix. The solution of this equation is of the general form
$\mathfrak{X}=\mathfrak{T}\mathfrak{J}_x\mathfrak{T}^{-1}$ with a transformation matrix $\mathfrak{T}$
and a Jordan matrix $\mathfrak{J}_x$, which has only zeros in the main diagonal. (accordingly to \cite{ga}, chapter
8.7) So that $\mathfrak{X}$ is isomorphic to a quaternion then  $\mathfrak{J}_x$ is of type (\ref{J4}) or
(\ref{Jk}). In the first case $\mathfrak{X}$ is the Zero matrix or isomorphic with it $\mathbf{X}$ is the
Zero quaternion (the trivial solution).\\
In the second case we have
$\mathfrak{X}=\mathfrak{T}\mathfrak{J}_x\mathfrak{T}^{-1}=\mathfrak{T}\left(\begin{array}{*{3}{c@{\;}}c}
	        0 & 1 & 0 & 0\\
			0 &  0 & 0 &  0\\
			0 &  0 &  0 & 1\\
			0 & 0 &  0 &  0
           \end{array}\right)\mathfrak{T}^{-1}$.
If the vector $\mathbf{x}$ was computed, then we can compute a special matrix according to (\ref{trans21}),
(\ref{trans22}) or (\ref{trans23}). Hence the solution is 
\begin{equation}\label{Lösung00}\mathbf{X}=\{\quat{x}\mid
x_0=0, \quad \mathbf{x}^2=x_1^2+x_2^2+x_3^2=0\}\end{equation}
The trivial solution is included in this solution set. In contrast to the solution set at non-singular quaternions
we do not have ${n \choose 2}$ spheres  but only one.\\ 
2. $\mathbf{a}\neq\mathbf{0}$ and $\mathbf{a}^2\neq 0$. (Case (\ref{Jd}))\\
So that $\mathfrak{A}$ is singular, we have the condition $a_0^2 +\mathbf{a}^2=0$.
Then $\mathfrak{A}$ has the eigenvalues $2a_0$ and 0 and the Jordan decomposition is
$\mathfrak{A}=\mathfrak{U}\mathfrak{J}\mathfrak{U}^{-1}$ with
$\mathfrak{J}=\left(\begin{array}{*{3}{c@{\;}}c}
	        2a_0 & 0 & 0 &\: 0\\
			0 &  2a_0 & 0 &\: 0\\
			0 &  0 &  0 &\: 0\\
			0 & 0 &  0 &\:  0
      \end{array}\right)$
If we restrict to solution matrizes that are isomorphic to a quaternion,this leads to type (\ref{Jd}) and we get: 
$\mathbf{X}=\sqrt[n]{2a_0}\left(
	\begin{array}{c}\frac{1}{2}\\[1ex]
	\frac{\mathbf{a}}{2a_0}
	\end{array}\right)=\frac{1}{2}\sqrt[n]{2a_0}\left(
	\begin{array}{c}1\\[1ex]
	\frac{\mathbf{a}}{a_0}
	\end{array}\right)$.\\
Here we have $n$ different values for $\frac{1}{2}\sqrt[n]{2a_0}$.\\
Example: Let $\mathbf{A}=\left(\begin{array}{c}
                         -8+8i\\
						 -16+8i\\
						 -8+24i\\
						 24+16i
						 \end{array}\right)\quad n=3$
The solutions are: \[\mathbf{X}_1=\left(\begin{array}{c}
											1+i\\
											1+2i\\
											3+i\\
											2-3i\end{array}\right)\quad
			\mathbf{X}_2=\frac{1}{2}\left(\begin{array}{c}
											-1-\sqrt{3}+(-1+\sqrt{3})i\\
											-1-2\sqrt{3}+(-2+\sqrt{3})i\\
											-3-\sqrt{3}+(-1+3\sqrt{3})i\\
											-1+3\sqrt{3}+(3+2\sqrt{3})i\end{array}\right)\]	
			\[\mathbf{X}_3=\frac{1}{2}\left(\begin{array}{c}
											-1+\sqrt{3}+(-1-\sqrt{3})i\\
											-1+2\sqrt{3}+(-2-\sqrt{3})i\\
											-3+\sqrt{3}+(-1-3\sqrt{3})i\\
											-1-3\sqrt{3}+(3-2\sqrt{3})i\end{array}\right)\]										
3. $\mathbf{a}\neq\mathbf{0}$ and $\mathbf{a}^2= 0$. (Kategorie (\ref{Jk}))\\
So that $\mathfrak{A}$ is singular then is $a_0=0$. It follows the jordan normal form of $\mathfrak{A}$.
$\mathfrak{A}=\mathfrak{U}\mathfrak{J}\mathfrak{U}^{-1}$ mit
$\mathfrak{J}=\left(\begin{array}{*{3}{c@{\;}}c}
	        0 & 1 & 0 & 0\\
			0 &  0 & 0 & 0\\
			0 &  0 &  0 & 1\\
			0 & 0 &  0 &  0
 \end{array}\right)$ and $\mathfrak{U}$ according to (\ref{trans21}), (\ref{trans22}) or (\ref{trans23}).
$\mathfrak{A}$ is nilpotent of index $2$. The solution matrix $\mathfrak{X}$ has to be nilpotent of index $x$
with $n<x\leq 2n$. We have $n\geq 2$, therefore $x\geq 3$. It follows that $\mathfrak{X}$ cannot be a  quaternion 
isomorphic matrix.
The equation $\mathbf{X}^n=\mathbf{A}$ is insoluble.
Example: $\mathbf{A}=\left(\begin{array}{c}0\\ 3\\ 4\\ 5i\end{array}\right)$
\section{Summary}
Now the results are displaied shortly: 
The quaternion equation $\mathbf{X}^n=\mathbf{A}$,  $n\geq 2$, integer, has the the following solution sets:\\
1a. $\mathbf{A}=a_0\mathbf{E}=a_0\einsq$, $a_0\neq 0$ \\
a) $\mathbf{X}=\sqrt[n]{a_0}\mathbf{E}$, $n$ different values of $\sqrt[n]{a_0}$ are possible.\\
b) Let $w_1$ and $w_2$ be 2 different n-th roots of $a_0$ ($w_1^n=w_2^n=a_0 $ but $w_1\neq w_2$). 
In order to compute $\mathbf{X}$ we have to consider all combinations of $w_1$ and $w_2$:\\
$\mathbf{X}=\{\vekk{x}\mid x_0=\frac{1}{2}(w_1+w_2)$, $\mathbf{x}^2=x_1^2+x_2^2+x_3^2=-\frac{1}{4}(w_1-w_2)^2\}$;
${n\choose 2}$ Kombinationen\\
1b. $\:\mathbf{A}=a_0\mathbf{E}=a_0\einsq$, $a_0= 0$ \\
$x_0=0, \quad \mathbf{x}^2=x_1^2+x_2^2+x_3^2=0$\\
2a. $\mathbf{A}=\vekk{a}$, mit $\mathbf{a}\neq\mathbf{0}$, $\mathbf{a}^2\neq 0$ und $\mathbf{a}^2+a_0^2\neq 0$.
\[\mathbf{X}=\sqrt[n]{a_0+i\sqrt{\mathbf{a}^2}}\left(
	\begin{array}{c}\frac{1}{2}\\[1ex]
	\frac{1}{2i}\frac{\mathbf{a}}{\sqrt{\mathbf{a}^2}}
	\end{array}\right)+\sqrt[n]{a_0-i\sqrt{\mathbf{a}^2}}\left(
	\begin{array}{c}\frac{1}{2}\\[1ex]
	-\frac{1}{2i}\frac{\mathbf{a}}{\sqrt{\mathbf{a}^2}}
	\end{array}\right)\]
According to the multiplicity of the n-th root there are $n^2$ solutions to $\mathbf{X}$.\\
2b. $\mathbf{A}=\vekk{a}$, with $\mathbf{a}\neq\mathbf{0}$, $\mathbf{a}^2\neq 0$
and $\mathbf{a}^2+a_0^2= 0$ (singular case)\\
$\mathbf{X}=\sqrt[n]{2a_0}\left(
	\begin{array}{c}\frac{1}{2}\\[1ex]
	\frac{\mathbf{a}}{2a_0}
	\end{array}\right)=\frac{1}{2}\sqrt[n]{2a_0}\left(
	\begin{array}{c}1\\[1ex]
	\frac{\mathbf{a}}{a_0}
	\end{array}\right)$ \\
According to the multiplicity of the n-th root there are $n$ solutions to $\mathbf{X}$.\\
3a. $\mathbf{A}=\vekk{a}$, with $\mathbf{a}\neq\mathbf{0}$, $\mathbf{a}^2= 0$ and $a_0\neq 0$\\
$\mathbf{X}=\sqrt[n]{a_0}\einsqv +\frac{1}{n\sqrt[n]{a_0^{n-1}}}\vekd{a}$\\
According to the multiplicity of the n-th root there are $n$ solutions to $\mathbf{X}$.\\
3b. $\mathbf{A}=\vekk{a}$, with $\mathbf{a}\neq\mathbf{0}$, $\mathbf{a}^2= 0$ and $a_0= 0$ (singular case)\\
In this case we do not have a solution.
			  
\end{document}